\documentclass[12pt,fleqn]{article}
\usepackage{color}
\usepackage{epsfig}
\usepackage{amsmath}
\usepackage{amssymb}
\usepackage{graphicx}
\usepackage[german,spanish,english]{babel}

\newcommand{\onehalf}{\mbox{$\frac{\scriptstyle 1}{\scriptstyle 2}\,$}}
\newcommand{\R}{{\mathbb R}}
\newcommand{\C}{{\mathbb C}}

\newtheorem{theorem}{Theorem}[section]
\newtheorem{lemma}[theorem]{Lemma}
\newtheorem{proposition}[theorem]{Proposition}
\newtheorem{corollary}[theorem]{Corollary}
\newtheorem{definition}[theorem]{Definition}
\newtheorem{example}[theorem]{Example}

\newenvironment{Theorem}{\begin{theorem} \begin{sl}}{\end{sl}
                          \end{theorem}}
\newenvironment{Lemma}{\begin{lemma} \begin{sl}}{\end{sl} \end{lemma}}
\newenvironment{Proposition}{\begin{proposition}
         \begin{sl}}{\end{sl} \end{proposition}}
\newenvironment{Corollary}{\begin{corollary}
       \begin{sl}}{\end{sl} \end{corollary}}

\newenvironment{Example}{\begin{example}
        \begin{rm}}{\end{rm} \end{example}}

\begin{document}

\title{Reduction by invariants, stratifications, foliations, fibrations and relative equilibria, a short survey.}

\renewcommand{\baselinestretch}{1}
\author { {\protect\normalsize J.C. van der Meer}\\ {\small\sl Faculteit
Wiskunde en Informatica, Technische Universiteit  Eindhoven},\\ {\small\sl PObox 513, 5600 MB Eindhoven, The
Netherlands}.}
\maketitle

\begin{abstract}
\noindent
In this note we will consider reduction techniques for Hamiltonian systems that are invariant under the action of a compact Lie group $G$ acting by symplectic diffeomorphisms, and the related work on stability of relative equilibria. We will focus on reduction by invariants in which case it is possible to describe a reduced phase space within the orbit space by constructing an orbit map using a Hilbert basis of invariants for  the symmetry group $G$. Results considering the stratification, foliation and fibration of the phase space and the orbit space are considered. Finally some remarks are made concerning relative equilibria and bifurcations of periodic solutions.
We will combine results from a wide variety of papers. We obtain that for the orbit space the orbit type stratification coincides with the Thom-Boardman stratification and each stratum is foliated with symplectic leaves. Furthermore the orbit space is fibred into reduced phase spaces and each reduced phase space has a orbit type stratification, where each orbit type stratum is a symplectic leaf.
\end{abstract}

\section{Introduction}
In Hamiltonian systems the symmetry is usually given by a  symmetry group acting on the phase space and leaving the system invariant. The symmetry group determines the geometry of the phase space and therefore also determines the behavior of systems having this symmetry. One of the tools to understand symmetric systems is reduction by dividing out the symmetry, which allows to study the dynamics on a lower dimensional space. The main questions are then how to construct this reduced phase space, how this reduced dynamics reconstructs to the original unreduced phase space, and which part of the dynamics persists under non-symmetric perturbations.

In this note we will consider reduction, including singular reduction. After a short review we will focus on reduction by invariants. By constructing an orbit map also singular reduction can be considered. Given a symmetry group $G$, fulfilling the right conditions, one can consider the orbit type stratification of the phase space and orbit space, the foliation by symplectic leaves of the orbit space, and the fibration of the phase space  by group orbits. Given a Hamiltonian system with symmetry group $G$ one can also consider relative equilibria. In general a relative equilibrium is an equilibrium position in a moving system. In the case of a symmetric system a relative equilibrium is an equilibrium that moves with the group action into a solution of the system and therefore corresponds to a stationary point of the reduced system. Relative equilibria for symmetric systems are of importance because, using singularity theoretic methods, it can be show that, under certain conditions, they persist under non-symmetric perturbation. As they are usually organised in families, they form an organising skeleton in the phase space of the symmetric system and perturbations thereof.

In this note we will focus on reviewing the different concepts of reduction, the stratification, foliation and fibration of the phase space and orbit space, and on the concepts of relative equilibria and their stability. We will give examples throughout the text and some more in the final section.

\section{Reduction}

The idea of reduction seems to go back on Reeb \cite{reeb1951,reeb1952} who studied perturbations of systems of differential equations only having periodic solutions. He considered the phase space as a fibred manifold or more precisely as an $S^1$ fibre bundle. He constructed the reduced system by mapping to the base space of the fibre bundle. One could call this concept fibre bundle reduction. He considered Hamiltonian systems as a special case.

Later the geometric reduction of Meyer \cite{meyer1973} and Marsden and Weinstein \cite{marsden1974} became the fundamental reference for reduction. Ideas on reduction can also be found in \cite{arnold1966} and \cite{smale1970}. The Meyer-Marsden-Weinstein reduction is formulated in \cite{marsden1974} as follows.

Consider a symplectic manifold M, on which we have the action by symplectic diffeomorphisms of a Lie group $G$ with Lie algebra $\mathfrak{g}$. For $\xi \in \mathfrak{g}$ let $\xi_M$ denote the corresponding infinitesimal generator or vector field on $M$. Let ${\cal J}:M\rightarrow \mathfrak{g}^*$ be an $Ad^*$ equivariant momentum map, that is, ${\cal J}\circ \varphi_g= (Ad_{g^{-1}})^*\circ {\cal J}$, where $\varphi_g$ is the action of $g\in G$, and $(Ad_{g^{-1}})^*$ is the co-ad{\cal J}oint action of $G$ on $\mathfrak{g}$. Let $\mu \in \mathfrak{g}^*$ be a regular value of ${\cal J}$, and let $G_{\mu}$ be the isotropy subgroup of $\mu$ for the co-adjoint action, and assume that $G_{\mu}$ acts properly and freely on ${\cal J}^{-1}(\mu)$. Then $M_{\mu}={\cal J}^{-1}(\mu)/G_{\mu}$ is the reduced phase space. If $\omega$ is the symplectic structure on $M$ then there exists an unique symplectic structure $\omega_{\mu}$ on $M_{\mu}$ with $\pi^*_{\mu}\omega_{\mu}=i^*_{\mu}\omega$, with $i^*_{\mu}$ the inclusion map of ${\cal J}^{-1}(\mu)$ in $M$, and $\pi^*_{\mu}$ the projection map of ${\cal J}^{-1}(\mu)$ onto $M_{\mu}$. If $H$ is a $G$-invariant Hamiltonian function on $M$ with respect to $\omega$, then the flow of $H$ induces a flow on $M_{\mu}$ which is Hamiltonian with respect to the symplectic form $\omega_{\mu}$ with reduced Hamiltonian function $H_{\mu}$.

Although this is a very general description it does not really allow us to construct the reduced phase space. Methods to construct the reduced phase space were given by Cushman \cite{cushman1981}, Van der Meer \cite{vandermeer1983, vandermeer1985}, Churchill, Kummer and Rod \cite{churchill1983}, Cushman an Rod \cite{cushman1982}, Kummer \cite{kummer1986}. They all use invariants for the group action to construct the reduced phase space. This method is formulated in a more general way using orbit spaces in \cite{abud1983,vandermeer1985} were a theorem of Hilbert is used that provides a Hilbert basis of invariants for compact group actions. This makes it possible to define an orbit map $\rho$ \cite{poenaru1976} under which all the $G$-orbits are mapped to points. Then $\rho({\cal J}^{-1}(\mu))$ is the reduced phase space, hence the name orbit space reduction or reduction by invariants. By a theorem of Schwarz $G$-invariant functions correspond to functions on the orbit space, that is, any $G$-invariant Hamiltonian function reduces to a function on the orbit space that naturally restricts to the reduced phase space. This method has the advantage that $\mu$ need not be a regular point, and that the action need not be free. The singularities of the orbit map reflect the fixed points, isotropy subgroups and orbit types. Note that singular reduction appears in \cite{vandermeer1983, vandermeer1985}
 for the non-semisimple 1:-1 resonance and in \cite{kummer1986} for the k:l resonances. This was later formalized in the context of Marsden-Weinstein reduction in \cite{arms1991}. We will distinguish between the Meyer-Marsden-Weinstein reduction and the construction using orbit maps by calling the first momentum map reduction and the latter reduction by invariants. Note that when using reduction by invariants the reduced phase space need not be a manifold. In general it is defined as a semi-algebraic set by relations and inequalities for the invariants defining the orbit map. Cushman and Sniaticky studied these spaces in more detail including differential structures on them (see \cite{cushman2001, cushman2005, sniatycki2013}) which allows to study more general vector fields on orbit spaces \cite{bates2021}.

The most general formulation for orbit map reduction does not start with a symplectic manifold but with a Poisson manifold. we consider $C^\infty(M)$ together with a Poisson bracket $\{\, ,\,\}$ making  $C^\infty(M)$ into a Lie algebra and call $(M, \{\, ,\,\})$ a Poisson manifold. If the Poisson structure is non-degenerate then $M$ is a symplectic manifold. The Poisson structure on $C^\infty(M)$ induces a Poisson structure on the orbit space.

Kummer \cite{kummer1981} gives a construction of the reduced phase space using principle G-bundles. In a neighborhood of regular values of the orbit map the reduced phase space is locally a $G$-bundle. That is, his ideas apply to $\rho({\cal J}^{-1}(\mu))$ without its critical set. The bundle structure of $G$-spaces is extensively studied in \cite{duistermaat2000}.

Reduction by invariants through the orbit map should not be confused with orbit reduction as introduced by Ortega \cite{ortega2002,marsden2007}. Orbit reduction refers to the fact that the pre-image of the co-adjoint orbit in the image of the momentum map is reduced. Reduction by invariants should also not be confused with Poisson reduction. Although Poisson reduction considers reduction of Poisson manifolds it restricts to $M=T^*G$ and its momentum map.

In orbit map reduction the orbit map is used to construct the reduced phase space and the reduced system. When using orbit map reduction one of the difficulties is to determine the invariants and the relations between the invariants. An extensive study of this is made in \cite{gatermann2000} (see also \cite{rumberger2001,chossat2002a}).

\section{Meyer-Marsden-Weinstein reduction, reduction by invariants and dual pairs}

Throughout this paper we will consider a connected, compact Lie group acting smoothly and properly on $\R^n$ which is assumed to be a symplectic space with the standard Poisson structure, that is it is a Poisson manifold, and $G$ is assumed to act by Poisson (symplectic) diffeomorphisms. Thus we assume to be in the nicest possible situation where we can use the strongest possible results. However, in many applications this is the case. Many of the results stated below will also hold under weaker conditions, more generally when $\R^n$ is replaced by a connected compact Poisson manifold $M$.

We will start by introducing reduction by invariants. In \cite{hilbert1890} Hilbert showed that
the algebra of polynomials over $\C$ of degree $d$ in $n$ variables which are invariant under $GL(n,\C )$, acting by substitution of variables, is finitely generated. This was extended by Weyl in \cite{weyl1946} who proved that
the algebra of invariants is finitely generated for any representation of a compact Lie group or a complex semi-simple Lie group.

Let $\R [x]_G$ denote the space of $G$-invariant polynomials with  coefficients in $\R$.
Consider a compact Lie group $G$ acting linearly on $\R^n$. Then there exist finitely many polynomials $\rho_1 ,\cdots ,\rho_k\in \R [x]_G$ which generate $\R [x]_G$ as an $\R$ algebra. These generators can be chosen to be homogeneous of degree greater then zero.
We call $\rho_1 ,\cdots ,\rho_k$ a Hilbert basis for $\R [x]_G$.

Schwarz \cite{schwarz1975} proved that if $\rho_1 ,\cdots ,\rho_k$ is a Hilbert basis for $\R [x]_G$, and
$\rho:\R^n\rightarrow \R^k; x\rightarrow (\rho_1(x) ,\cdots ,\rho_k(x))$. Then $\rho ^*: C^\infty (\R^k,\R ) \rightarrow C^\infty(\R^n,\R )_G$ is surjective, with $\rho^*$ the pull-back of $\rho$.
Thus all $G$-invariant smooth function can be written as smooth functions in the invariants.

The following, showing that $\rho$ is an orbit map, can be found in \cite{poenaru1976}.
The map $\rho$ is proper and separates the orbits of $G$. Moreover the following diagram commutes, with $\tilde{\rho }$ a homomorphism
\begin{alignat*}{2}
\R^n & \;\;\stackrel{\rho}\longrightarrow &\rho(R^n)\\
\pi\searrow & &\!\!\!\swarrow\tilde{\rho } \\
        &\; \; \R^n/G &
\end{alignat*}

Here the orbit space $\R^n/G$ is the quotient space $\R^n/ \sim$, where the equivalence relation is given by $x\sim y$ if $x$ and $y$ are in the same $G$-orbit. We can take $\rho(\R^n)$ as a model for the orbit space.

Consider $(\R^{2n},\omega)$ on which a Lie group $G$ acts linearly and symplectically. Then $(C^\infty (\R^{2n},\R ), \{\; ,\;\})$ is a Poisson algebra. If we consider on $\R^k$ the Poisson structure induced by $\rho$ by taking as structure matrix $W_{ij}=\{\rho_i,\rho_j\}$ then
$(C^\infty (\R^k,\R ), \{\; ,\;\}_W)$ is a Poisson algebra and $\rho$ a Poisson map. We have a reduction of the Poisson manifold if we restrict the bracket on $\R^k$ to $\rho (\R^{2n})$.

In general there will be relations and inequalities determining the image of $\rho$. Therefore $\rho (\R^{2n})$ will in general be a real semi-algebraic subset of $\R^{k}$, where a semi-algebraic subset of $\R^k$ is a finite union of sets of the form
$\{x\in \R^k|R_1(x)= \cdots =R_s(x)=0\; , \; R_{s+1}(x),\cdots ,R_m(x) \geqslant0\}$ .
Define $C^\infty (\rho (\R^{2n}),\R)=\{ F:\rho (\R^{2n})\rightarrow \R|\rho^* (F)\in C^\infty(\R^{2n},\R)\}$. This is a differential structure on $\rho (\R^{2n})$ and  the orbit map is smooth (see \cite{cushman2001, cushman2005}). Note that the $R_i$, $1\leqslant i\leqslant s$, are Casimirs for the induced Poisson structure $\{\; ,\;\}_W$.

Let $W$ be a real semi-algebraic variety in $\R^k$.
A point $x\in W$ is nonsingular if there exists a neighborhood $U\subset W$ of $x$ such that for each $y\in U$ the matrix
$\frac{\partial R_i}{\partial x_j}(x)$ has maximal rank.
A point $x\in W$ is singular if the rank of $\frac{\partial R_i}{\partial x_j}(x)$ is strictly less than the maximal rank.

Combining remarks in \cite{marsden2007} and \cite{simo1991}  we find that for any Poisson manifold $(M, \{\, ,\,\})$ on which we have a compact Lie group $G$ acting by Poisson maps, $G$ has a Lie algebra $\mathfrak{g}$. To each element $\xi \in \mathfrak{g}$ we may associate a Hamiltonian vector field $X_{J_\xi}=\{J_\xi , \cdot \}$ on $M$ with Hamiltonian function $J_\xi$ defined by $J(z)={\cal J}(z)\cdot \xi$.

Let $\xi_i$ be a basis of $\mathfrak{g}$ such that $exp(\xi_i)$ generate $G$. Consider the corresponding functions $J_i$ and consider the momentum map $J(z)=(J_1(z), \cdots ,J_r (z))$.
Furthermore consider the Hilbert basis of invariants $\rho_i, \; i=1, \cdots ,k,$ for the $G$-action. Obviously, $\{ \rho_i,J _j\}=0$ for all $i,j$. Thus the maps $\cal{J}$ and $\rho$ form a dual pair. Note that without further conditions the images of $\cal{J}$ and $\rho$ are at best semi-algebraic sets. Also the $\rho_i$ need not span a Lie algebra. however, the $\rho_i$ generate the Poisson algebra of $G$-invariant functions $C^\infty (M)^G$. Let $G'$ denote the Lie group generated by the G-equivariant Poisson vector fields $X_f$, $f \in C^\infty (M)^G$. Then the map $M\rightarrow M/G'$ is called the optimal momentum map\cite{ortega2002,ortega2004}. If the $J_i$ are the invariants defining  $M/G'$ the $J$ is optimal. Consider $C^\infty (J(M))$, then $C^\infty (M)^G$ and $C^\infty (J(M))$ centralize each other in the Poisson algebra $C^{\infty} (M)$ (see for instance \cite{karshon1997}), i.e. we have a Howe dual pair \cite{howe1985}. Actually, as can be found in \cite{ortega2003}, the pair $\mathfrak{g}\leftarrow M \rightarrow \rho (M) $ is a Lie-Weinstein dual pair.

As $\{ \rho_i,J _j\}=0$ it follows that

\begin{Proposition}\label{kerdJ}
$ker\,(\text{d}J)$ is spanned by the Hamiltonian vector fields $X_{\rho_i}$.
\end{Proposition}

\begin{Example}\label{symplso3}
Consider $SO(3)$ acting on $\R^3$. The lifted action on the co-tangent bundle $T^*\R^3$ is the diagonal action of $SO(3)$ on $\R^6$. Let $(x,y)$ denote the coordinates on $\R^6=\R^3\times\R^3$. The generators for the group are $x_1y_2-x_2y_1, x_1y_3-x_3y_1, x_2y_3-x_3y_2$ (or the components of the cross product $x\times y$), which span a Lie algebra isomorphic to $\mathfrak{so}(3)$. This action is a symplectic action on $\R^6$ with respect to the standard symplectic form. A Hilbert basis for this action is $|x\times y|^2$, and $|x|^2, |y|^2, <x,y>$. The last three invariants span a Lie algebra isomorphic to $\mathfrak{sl}(2,\R)$. According to \cite{weyl1946} in general the full linear Lie algebra invariant under a diagonal $SO(n)$-action is $\mathfrak{s}l(2,\R)$. Consequently the momentum map
\[
J:(x,y)\rightarrow (x_1y_2-x_2y_1, x_1y_3-x_3y_1, x_2y_3-x_3y_2)\; ,
\]
and the orbit map $\sigma:(x,y)\rightarrow (|x|^2, |y|^2, <x,y>)$ are a dual pair. The orbit space is defined by Lagrange identity $ |x|^2|y|^2-<x,y>^2=|x\times y|^2$ together with the inequalities $|x|^2\geqslant 0$, and $|y|^2\geqslant 0$. It is a solid cone. The reduced phase spaces are given by taking $x\times y$ constant, thus, provided $|x\times y|\neq 0$, the reduced phase space is one sheet of a two sheeted hyperboloid. In case $|x\times y|= 0$ it is a cone.
\end{Example}

\begin{Example}
Consider an integrable system on $\R^{2n}$, with $n$ integrals in involution, that is, the group $G$ generated by the integrals is a torus.
The momentum map and the orbit map are the same, i.e. $J_i=\rho_i$. The orbit space is a polytope \cite{atiyah1982, guillemin1984}. A reduced phase space is a point. Regular points correspond to $n$-tori. The faces, edges and vertices of the polytope are the images of the singular points of the orbit map and correspond to lower dimensional tori.
\end{Example}

\begin{Example}
Consider a group $G$ which is the flow of a linear Hamiltonian system in two degrees of freedom in $k:\ell$ resonance, $k\in \mathbb{N}$, $\ell\in \mathbb{Z}$, $|k|\neq \ell$). That is, the Hamiltonian generating the group $G$ is $H(x,y)=\onehalf k(x_1^2+y_1^2)+\onehalf \ell (x_2^2+y_2^2)$. The action of $G$ is an $S^1$-action and a Hilbert basis for the invariants of $G$ is given by polynomials $I_1(x,y)$, $I_2(x,y)$, $R_1(x,y)$ and $R_2(x,y)$, where $I_1,I_2$ are quadratic polynomials and  $R_1,R_2$ are polynomials of order $k+|\ell|$. When we choose $I_1(x,y)=H(x,y)$ and $I_2(x,y)=\onehalf k(x_1^2+y_1^2)-\onehalf l(x_2^2+y_2^2)$ we have bracket relations $\{ I_2,R_1\}=2klR_2$, $\{ I_2,R_2\}=2klR_1$, $\{ I_1,R_1\} =\{ I_1,R_2\} =0$, and $\{ R_1, R_2\} = -2(I_1^2-I_2^2)+(I_1+I_2)^2$. Furthermore we have the following relation $R_1^2+R_2^2=\onehalf (I_1+I_2)^2(I_1-I_2)$ (see \cite{sanders1992}). In this case the invariants do not form a Lie subalgebra of $C^\infty (\R^4)^G$. However, we do have an orbit map $\rho:(x,y)\rightarrow (I_1,I_2,R_1,R_2)$, where the image is determined by $R_1^2+R_2^2=\onehalf (I_1+I_2)^2(I_1-I_2)$, $I_1\geqslant 0$. The reduced phase space is obtained by taking $I_1(x,y)=c$, giving the reduced phase space given by the relation $R_1^2+R_2^2=\onehalf (c+I_2)^2(c-I_2)$ in $(I_2,R_1,R_2)$-space (cf \cite{kummer1986}).
\end{Example}

We have a special situation if the Hilbert basis is a finite Lie subalgebra of $C^\infty (\R^n)$. In this case the orbit map can also be interpreted as a momentum map. We have a Lie-Weinstein dual pair of momentum maps.
\[
\mathfrak{g}_1^*\leftarrow M \rightarrow \mathfrak{g}_2^* \; ,
\]
or in terms of generating functions
\[
(J_1, \cdots ,J_r)\leftarrow M \rightarrow (\rho_1, \cdots ,\rho_k) \; .
\]
If the $J_i$ generate $G_1$ and the $\rho_i$ generate $G_2$ then $G_1$ and $G_2$ commute. Moreover $C^\infty (M)_{G_2}$ is the Lie algebra of smooth functions in the $J_i$.
The center of $\mathfrak{g}_1$ therefore consists of functions in $C^\infty (M)_{G_2}$.

This holds more generally. If we have a group $G$, with momentum map $J$, then the center of $C^\infty (M)^G$ in the Poisson algebra $C^\infty (M)$ consists of smooth functions in the $J_i$, which can be considered as the universal enveloping algebra of $\mathfrak{g}$. As a consequence the Casimirs of $C^\infty (M)_G$ do  belong to $C^\infty (M)_G$ and to the universal enveloping algebra of $\mathfrak{g}$.
Thus the symplectic leaves for the Poisson structure on the image of the orbit map, which are obtained by setting the Casimirs equal to a constant \cite{weinstein1983}, are given by $\rho(J^{-1}(\mu))$.

When the invariants form a Lie algebra then the symplectic leaves are the co-adjoint orbits of $G_2$ on $\mathfrak{g}_2$ \cite{weinstein1983}.

\section{Stratifications, foliations and fibrations}
Orbit spaces were studied in connection to understanding the structure of $G$-spaces in the 1950's \cite{palais1960}. Mainly to understand extrema of $G$-invariant functions \cite{michel1971}, \cite{ michel1972},\cite{abud1983} in connection to applications in solid state physics.  More recent applications in quantum mechanics and molecular behavior that connect to the subject of this paper are \cite{michel2001}, \cite{zhilinskii2001}, \cite{sadovskii2018}. In \cite{watts2013} a very nice overview of all the relevant theorems from the literature is given in connection to orbit spaces in the context of Meyer-Marsden-Weinstein reduction. When studying symmetric spaces and reduction one of the important issues is the orbit type stratification, see the above cited literature and \cite{sjamaar1991}, \cite{duistermaat2000}, \cite{audin2004} and the very accessible notes \cite{duistermaat2008}, \cite{meinrenken2003}.

Consider a Hamiltonian $G$-action on a connected manifold $M$, $G$ a compact Lie group acting properly and smoothly on $M$.
Let $G\cdot x=\{y\in M|y=g\cdot x, g\in G\}$ be the $G$-orbit in $M$ through $x\in M$.
$G_x=\{g\in G|g\cdot x=x$ is the isotropy subgroup of $G$ at $x$. $G_x$ is a closed Lie subgroup of $G$. If $H$ is a subgroup of $G$ for which $H_x=\{g\in H|g\cdot x=x\}$, then we may call $H_x$ also an isotropy subgroup of $x$. Obviously $H_x \subset G_x$ and $G_x$ is the maximal isotropy subgroup for $x$ in $G$. We have
\[
\textrm{dim} (G)=\textrm{dim}(G_x) +\textrm{dim} (G\cdot x) \; .
\]
We have $G_{g\cdot x}=gG_x g^{-1}$, that is, the isotropy subgroups of points on the same orbit are each others conjugate and therefore isomorphic.
For a subgroup $H$ of $G$ define the normalizer of $H$ by
\[
N_G(H)=\{g\in G|gHg^{-1}=H\} \; .
\]
$N_G(H)$ is a closed Lie subgroup of $G$ and the largest subgroup of $G$ containing $H$ as a normal subgroup.

Let $H$ be a subgroup of $G$. Let $M_H=\{x\in M|G_x=H\}$. $M_H$ is called the isotropy type of $H$ in $M$.

\begin{Lemma} \cite{duistermaat2008}\label{duislem1}
Let $H$ be a closed Lie subgroup of $G$. Then the action of $N(H)$ leaves $M_H$ invariant and induces a free action of $N(H)/H$ on $M_H$.
\end{Lemma}

Denote by $[H]$ the conjugacy class of $H$. Set $M_{[H]}=\{m\in M|G_m=gHg^{-1}\, ,\, g \in G\}$. $M_{[H]}$ is called the orbit type of $[H]$ in $M$. Two points in $M$ belong to the same orbit type if and only if their exists a $G$-equivariant bijection between their $G$-orbits \cite{duistermaat2000}. $M$ is partitioned into orbit types $M_{[H]}$ and each orbit type is partitioned into isotropy types $M_K$, $K\in [H]$ for conjugate subgroups of $G$. This induces a partition of $\rho (M)$ into orbit types $\rho(M_{[H]})$.

\begin{Lemma}\cite{duistermaat2008}\label{duislem2}
$\rho (M_H)=\rho (M_{[H]})$ and because the $G$-action is proper we have that $\rho|_{M_H}:M_H\rightarrow \rho(M_{[H]})$ is a principal fiber bundle with structure group $N(H)/H$.
\end{Lemma}

Note that if $M$ is $\R^n$ and we have a faithful representation of $G$ on $\R^n$, with a subgroup $H$, then $M_H= Fix (H)$. On $Fix(H)$ the action of $G$ reduces to the action of $N(H)/H$.

It is important to note that $M_H$ can have connected components with different dimensions. If the action is free then the orbit space $M/G$ is a smooth manifold, there is only one orbit type, and the orbit map is a smooth fibration with structure group $G$. $\textrm{dim}(M/G)=\textrm{dim}(M)-\textrm{dim}(G)$

For instance in \cite{duistermaat2000} it can be found that the connected components of the orbit types form a Whitney stratification in $M$.

There exists a partial ordering of isotropy and orbit types. For isotropy subgroups $H$ and $K$ we say that $M_H\leqslant M_K$, $M_{[H]}\leqslant M_{[K]}$ if and only if $H$ is conjugate to a subgroup of $K$.

\begin{Theorem}[Principal orbit theorem]
As before consider a group action by $G$ on a connected differentiable manifold $M$. Then there exists a maximal orbit type. The maximal orbit type stratum $S_m$ is open and dense in $M$. The orbit space $S_m/G$ is open and dense and connected in $M/G$.
\end{Theorem}
The maximal orbit and its orbit type are also called the principal orbit and principal orbit type.

Denote the orbit type strata by $S_i$, that is, $M=\cup_i S_i$. When we have defined an orbit map by invariants we obtain an orbit type stratification of the orbit space $\rho(M)=\cup_i \rho(S_i)$.

Note that in general our manifolds are embedded in $\R^n$ and thus separable. In \cite{palais1960} we find the following

\begin{Corollary}\label{dim1}\text{}
\begin{itemize}
\item[(i)] $\rm{dim}(\rho (M_{[H]}))=\rm{dim}(M_{[H]})-\rm{dim}G/H \; , \; H\subseteq G$ ,\\
\item[(ii)] $\rm{dim}(M/G)=sup\{\rm{dim}(\rho(M_{[H]}))-\rm{dim}(G/H)|H\subseteq G\}$ .
\end{itemize}
\end{Corollary}

On the orbit space the partitioning given by the stratification consists of sets of points that have as a pre-image diffeomorphic orbits.
A G-orbit is the pre-image $\rho^{-1}(\nu)$ of a point $\nu$ on the orbit space. More precisely,  if $\nu \in \rho (M_H)$ the orbit through $p\in \rho^{-1}(\nu)$ is $N(H)/H \cdot p$. Its tangent space at a point $p\in \rho^{-1}(\nu)$ is spanned by the Hamiltonian vectors $X_{J_i}(p)$. We have
\[
\rm{dim}(G.p)=\rm{rank}(J)(p)= \rm{dim}(N(H)/H \cdot p)\; .
\]
The orbit space is given, as a subset of $\R^k$, by the relations and inequalities for the invariants. If the action of $G$ is free the principal orbit type corresponds to the isotropy subgroup $I$ and we have
by taking $H=I$ in corollary \ref{dim1}, that the dim($S_m)=$dim$(\rho (M))$,
\[
\textrm{dim} (G) +\textrm{dim} (\rho(M))= \textrm{dim} (M) \; ,
\]
where dim$(\rho(M))$ equals the maximal rank of $\rho$. Moreover, the dimension of the orbit type stratum $\rho(S_i)$ for some $G$-orbit in $M$ equals the rank of the orbit map at the points of this orbit. Consequently the orbit type stratification coincides withe the Thom-Boardman stratification \cite{golubitsky1973} for the orbit map.The image of this stratification under the orbit map is of course the same as that of the orbit type stratification, but now it can be seen as the singular set stratification of the semi-algebraic set that is te image of the orbit map. The stratification of the image of the orbit map defines a partition of the semi-algebraic set into disjoint sets on which the rank of $\rho$ is constant and on which each point corresponds to an orbit of the same type.

Again let $p\in \rho^{-1}(\nu)$. And let $\nu$ be in the orbit type stratum for $M_{H}$. Following from lemma's \ref{duislem1} and \ref{duislem2} we have that, $p\in M_H$, $\rm{dim}(N(H)/H)=\rm{dim}(G\cdot p)=\rm{rank}(J)(p)$. Furthermore
\[
\rm{dim}(G)=\rm{dim}(G\cdot p)+\rm{dim}(G_p) \; ,
\]
and
\[
\rm{dim}(M_H)=\rm{dim}(G\cdot p)+\rm{rank}(\rho)(\nu) \; .
\]

So far we have paid attention to the reconstruction of the orbit space in terms of $G$-orbits and orbit type strata. The orbit space is fibred into reduced phase spaces $\rho (J^{-1}(\mu ))$. Each reduced phase space has a stratification into orbit type strata by considering the intersection of the reduced phase space with the orbit type strata, i.e. $\rho (J^{-1}(\mu )) \cap \rho (M_{[H]})$. This is the symplectic stratification introduced in \cite{sjamaar1991}.

On the other hand the Poisson structure on the orbit space allows us to obtain a foliation into symplectic leaves \cite{weinstein1983} of the orbit type strata. Each symplectic leaf is obtained by setting the Casimirs of the Poisson structure equal to a constant. Consequently, the symplectic leaves are subspaces of the reduced phase spaces. Besides that, the rank of the Poisson structure is constant along a symplectic leaf.

\begin{Proposition}
If $\tilde{W}$ is the invertible structure matrix of the Poisson structure on the symplectic manifold $M$, then $(d\rho )\tilde{W} (d\rho )^T$ is the induced structure matrix on the orbit space $\rho (M)$, and the rank of the induced Poisson structure on the orbit space equals the rank of the orbit map, that is,
\begin{equation}
\rm{rank} ((d\rho )\tilde{W} (d\rho )^T) = \rm{rank} (d\rho ) \; .
\end{equation}
\end{Proposition}

 The proof is a straightforward exercise in linear algebra. Now the stratification by rank of the orbit space coincides with the orbit type stratification, thus each orbit type stratum is foliated into symplectic manifolds of the same dimension on which the Poisson structure has the same rank. The same then also holds for the orbit type strata of the reduced phase spaces. Because the reduced phase space as well as the symplectic leaves are obtained by setting the Casimirs equal to a constant it follows that each orbit type stratum of the reduced phase space is a symplectic leaf.

Thus conclusively
\begin{Theorem}
On the orbit space the orbit type stratification coincides with the Thom-Boardman stratification and each stratum is foliated with symplectic leaves.
Furthermore the orbit space is fibred into reduced phase spaces and each reduced phase space has a orbit type stratification, where each orbit type stratum is a symplectic leaf for the induced Poisson structure.
\end{Theorem}
Recall that each orbit type stratum is a principal fibre bundle \ref{duislem2}.

When we consider a symplectic manifold $M$ which is embedded in $\R^n$, and which is a symplectic leaf for the Poisson structure on $\R^n$, then we may further reduce if a group $G$ is acting on  $\R^n$ by Poisson diffeomorphisms and leaving $M$ invariant. This way we may reduce in stages \cite{marsden2007}.

\section{Relative equilibria}
As before consider $\R^n$ with the standard non-degenerate Poisson structure (thus $n$ is even) and a group action of a compact and connected Lie group $G$ by symplectic (Poisson) diffeomorphisms. Let $H \in C^\infty (\R^n, \R )^G$ and consider the Hamiltonian (Poisson) vector field $X_H$ on $\R^n$. $X_H$ has integrals $J_i$. Because $H \in C^\infty (\R^n, \R )^G$ there is a function $\tilde{H} \in C^\infty (\R^k, \R )$  on the target space of the orbit map such that $H=\tilde{H} \circ \rho$. The reduced vector field is now the Poisson vector field $X_{\tilde{H}}$ with respect to the induced Poisson structure $\{ \, , \, \}_W$.

In \cite{abraham1978} we find as a definition for relative equilibrium that a point $x \in \R^n$ is a relative equilibrium for $X_H$ if $\rho(x)$ is a stationary point for the reduced vector field. Other ways of formulating this are that a relative equilibrium is a point $x \in \R^n$ such that the solution of Hamilton's equations for $X_H$ with initial value $x$ coincides with the orbit of a one parameter sub-group  of $G$ \cite{patrick1992} or, somewhat different, that relative equilibria are $G$ group orbits which are invariant under the flow of $X_H$ \cite{roberts2002}.

We will use the formulation from \cite{simo1991} stating that a point $x \in \R^n$ is a relative equilibrium for the Hamiltonian system $X_H$ if the trajectory $\gamma_t$ of Hamilton's equations for $X_H$ through $x$ is given by
\[
\gamma_t(x)= \textrm{exp}(t X_F )(x) \; ,
\]
where$X_F$ with $F=\sum_{i=1}^r \lambda_i J_i$ is an infinitesimal generator for an element of $G$ and $x \in J^{-1}(\mu)$.
Thus $X_H(x)=X_F (x)$. Thus a relative equilibrium is a critical point for the energy-momentum map $H\times J:\R^n\rightarrow \R^{r+1}; x\rightarrow (H(x),J_1(x),\cdots ,J_r(x))$ \cite{smale1970}.
In \cite{smale1970} a reduction is performed by considering $(H\times J)^{-1}(h,\mu)/G_x$.
Obviously, as a relative equilibrium is contained in a $G$-orbit it maps to a point $\rho (x)$ in the reduced phase space $\rho(J^{-1}(\mu))$ if $x\in J^{-1}(\mu)$. Furthermore as the $X_H$ trajectory through $x$ is a $G$-orbit it reduces to a stationary point $\rho(x)$ for $X_{\tilde{H}}$.

Of specific interest are of course the stability of relative equilibria and the persistence under change of the momentum. The latter results in families of relative equilibria that might be organized in manifolds. The study of these concepts has a long history. Some relevant references are \cite{simo1991}, \cite{patrick1992}, \cite{ortega1999}, \cite{montaldi1997}, \cite{lerman1998}, {\cite{patrick2004},  \cite{montaldi2011}. As for the results concerning stability up till that moment a very nice discussion is given in the introduction of \cite{patrick2004}.

In \cite{simo1991} the energy-momentum method is introduced to determine the formal stability of a relative equilibrium. To determine the relative equilibria we have to solve the Lagrange multiplier optimization problem of finding the critical points of $H$ under the constraints $J(x)=\mu$. To determine the stability one considers ${\rm d}^2H(x)$. However, one has to restrict ${\rm d}^2H(x)$ to some subspace $\cal{S}$ of $ker{\rm d}J(x)$. To determine this subspace one has to remove the neutral directions from $ker(\rm{d}J(x))$. By Lemma \ref{kerdJ} $ker({\rm d}J(x))$ is spanned by the $X_{\rho_i}(x)$. However, there might be dependencies at $x$. These dependencies determine the neutral directions and are given by the Casimirs that determine the symplectic leaf through $\rho (x)$. Note that these Casimirs can be expressed as smooth functions of the $\rho_i$ and as smooth functions of the $J_i$.  Suppose this set of Casimirs (independent at $\rho (x)$) is given by $C_i$, $i=1,\cdots, s$. Then the vectors $X_{C_i}(x)$ determine the directions to be left out of  $ker(\rm{d}J(x))$. The vector fields $X_{C_i}$ are the infinitesimal generators for the isotropy subgroup $G_\mu$ of $G$, which is the group leaving $\mu$ fixed under the co-adjoint action of $G$ on $\mathfrak{g}^*$. That is, ${\cal S} = ker({\rm{d}}J(x))/T_x (G_\mu\cdot x)$ as in \cite{simo1991}.

Now the vectors $X_{\rho_i}(x)$ reduce to tangent vectors to the symplectic leaf in the reduced phase space at $\rho (x)$. If one leaves out the dependencies given by the Casimirs we get that $\cal{S}$ maps to the tangent space to the symplectic leaf through $\rho(x)$ at $\rho (x)$. Thus
\begin{Corollary}
Formal stability of the relative equilibrium $x$ as defined in \cite{simo1991} agrees with stability of $\rho (x)$ on the reduced phase space.
\end{Corollary}
This is called relatively stable in \cite{abraham1978}.

The above becomes more clear when we consider the following to be found in \cite{gatermann2000}

\begin{Theorem} Consider a Hilbert basis $\pi_1,\cdots , \pi_r$ for a faithful representation of a compact Lie group. Assume $H$ is an isotropy subgroup of $G$ in this representation with corresponding Fix(H). Then there exist invariants $\tilde{\pi}_1, \cdots ,\tilde{\pi}_d$,
which are algebraically independent polynomials in the $\pi_i$, such that
\[
\tilde{\pi}_i|_{Fix(H)} \neq 0\; , i=1,\cdots ,c\; , \rm{and} \; \tilde{\pi}_i|_{Fix(H)}= 0\; , i=c+1,\cdots ,d\;'.
\]
furthermore $\tilde{\pi}_i|_{Fix(H)}$, $i=1,\cdots ,c$ are algebraically independent.
\end{Theorem}

That is, one can pick a set of invariants defining the orbit space such that the tangent space at a point $p$ of the reduced phase space is spanned by the tangent vectors generated by $\tilde{\pi}_i|_{Fix(H)} \neq 0\; , i=1,\cdots ,c$, where $H$ is the isotropy group of a point in $\rho^{-1}(p)$. (see also \cite{koenig1997}).

Thus the $\tilde{\pi}_i|_{Fix(H)} \neq 0\; , i=1,\cdots ,c$, can be considered as a set of invariants defining the orbit space for $Fix(H)$, while $\tilde{\pi}_i|_{Fix(H)} \neq 0\; , i=c+1,\cdots ,d$ can be considered as Casimirs, that is, the set of Casimirs becomes larger when there is a non-trivial isotropy subgroup, and consequently the orbit space and the reduced phase spaces reduce in dimension.

\begin{Example}
Consider an integrable system, i.e. $G$ is the torus group. The orbit map and the momentum map are the same thus the orbit space is the momentum polytope. The interior, faces, edges and vertices of this polytope correspond to the orbit type strata. The reduced phase spaces are points. Consequently each point is a relative equilibrium. As the pre-image of a point is a torus the trajectories of $G$-invariant vector fields are periodic orbits or quasi-periodic orbits. Quasi-periodic relative equilibria are considered in \cite{duistermaat2008}.
\end{Example}

In many examples the presentation of $G$ is linear and explicitly known. In \cite{montaldi1988} we find the possible representations that can occur when we consider the linear symplectic action of a Lie group $G$ on a vector space $V$.

\begin{Theorem}( \cite{montaldi1988}) Every symplectic representation $V$ of $G$ has a unique direct sum decomposition
\[
V=V_1 \oplus \cdots \oplus V_{\ell} \; ,
\]
where
\begin{itemize}
\item[(a)] The $V_j$ are $G$-invariant subspaces of $V$;
\item[(b)] $V_j= \mathbb{K}_j^{n_j}\otimes_{\mathbb{K}_j} W_j$, where $W_1,\cdots ,W_{\ell}$ are pairwise irreducible representations of $G$ and $Hom_G(W_j,W_j)\approx \mathbb{K}_j=\mathbb{R},\mathbb{C},\mathbb{H}$;
    \item[(c)] The action of $G$ on $V_j$ is the tensor product of the action $W_j$ and the trivial action on $\mathbb{K}_j^{n_j}$
\end{itemize}
\end{Theorem}

Let $Sp_G(\mathbb{R}^{2n})$ denote the group of $G$-equivariant symplectic linear transformations on $\mathbb{R}^{2n}$. Then \cite{montaldi1988} $Sp_G(\mathbb{R}^{2n})\cong S_1 \times \cdots S_{\ell}$, where each $S_j$ is either
$Sp(m,\mathbb{R})$, $U(p,q;\mathbb{C})$ or $\alpha U(r,\mathbb{H})$. Here $Sp(m,\mathbb{R})$, $U(p,q;\mathbb{C})$ and $\alpha U(r,\mathbb{H})$ are as defined in \cite{montaldi1988}. As we are dealing with linear symplectic maps the corresponding Lie algebra $\mathfrak{sp}_G(\mathbb{R}^{2n})$ is isomorphic to the Lie algebra under the Poisson bracket of $G$-invariant homogeneous quadratic polynomials.
If these polynomials form a Hilbert basis then we are in the situation of a dual pair $\mathfrak{g}$, $\mathfrak{sp}_G(\mathbb{R}^{2n})$.

\section{Bifurcations of periodic solutions}
When $G$ is a symplectic $S^1$ action the relative equilibria are periodic solutions. If furthermore the Hamiltonian depends on parameters one can study the bifurcation of these periodic solutions in dependence of the parameters. Here we have to distinguish between the parameters in the Hamiltonian, which are sometimes called distinguished parameters, or unfolding parameters, that usually are related to the physical system parameters, and parameters introduced by the value of the momentum map, i.e. introduced by the reduction. When considering a $G$-invariant system in the neighbourhood of a stationary point one can, if the quadratic part of the Hamiltonian fulfills certain conditions, consider the additional $S^1$ action by the semisimple part of this quadratic Hamiltonian. The bifurcation one wants to describe is then the bifurcation of periodic orbits with period close to the period of this $S^1$ action. To this end one first uses Liapunov-Schmidt reduction, or a splitting theorem, to reduce the $G$-invariant system to a $G\times S^1$-invariant system. Then $G\times S^1$-equivariant singularity theory is used to reduce the power series of the  $G\times S^1$-invariant Hamiltonian in the neighborhood of the stationary point to a finite part of the power series. Note that this depends on several non-degeneracy conditions that have to be fulfilled. When one finally has a $G\times S^1$ invariant polynomial system the parameter dependent equation for the stationary points of the $S^1$-reduced system then gives the bifurcation equation. For $S^1$-symmetric systems these ideas were introduced in \cite{duistermaat1984}, \cite{vandermeer1985}, and for $G\times S^1$ invariant systems in \cite{montaldi1988, montaldi1990a, montaldi1990b}. Also see \cite{golubitsky1985, golubitsky1988}. Note that on the $S^1$ bifurcation picture one still has the action of the group $G$, this leads to bifurcations with symmetry \cite{vandermeer1990}, \cite{chossat2002b}.

Following \cite{montaldi1988} we have for the representation of a linear symplectic $G$-action the decomposition $V=V_1 \oplus \cdots \oplus V_{\ell}$, i.e. $G=G_1 \times \cdots \times G_{\ell}$, where $G_i$ acts on $V_i$. For any component $G_i$ which is an $S^1$-action we may find periodic solutions on the fixed point spaces of the subgroups of $G_i$ using the ideas of \cite{montaldi1988}. Note that isomorphic subgroups might have different fixed point spaces in terms of geometric place. Thus for each $G_i$ which is an $S^1$-action one can pass to the reduced phase space for this $S^1$-action and find stationary points on all isotropy strata.

\section{Examples}

In many applications the symmetry group can be decomposed into $S^1$ (or $SO(2)$) symmetries and $SO(3)$ symmetries. In this section we will give some examples.

\subsection{Perturbed harmonic oscillator on $\R^8$}\label{harmosc}
In this example we will consider a perturbed harmonic oscillator on $\R^8$ modelling the Van der Waals and Zeeman system. In view of \cite{vandermeer2015}, \cite{vandermeer2021} these can be considered as models for a regularized perturbed Keplerian systems. These models can be found in \cite{diaz2010},\cite{egea2011} were the focus was on the bifurcation of relative equilibria.

Consider on $\R^8$ with standard symplectic form and coordinates $(q,Q)$ the Hamiltonian
\[
\begin{split}
\mathcal{H}(K,N,\Xi,L_{1}, H_2)=& H_2+\dfrac{3}{4}\left(  3\beta^2-2\right)
K^{2} H_2+ (1-\beta^2)K\Xi L_1+\dfrac{1}{2}\left(  4-\beta^2\right)  N H_2\\
&+(\dfrac{3}{2}+\dfrac{\beta^2}{4})H_2^3-  (\dfrac{\beta^2 }{2}+1)\dfrac{H_2}{2}\left(  {L_1}^{2}+{\Xi}^{2} \right) \; ,
\end{split}
\]
with
\begin{align*}
{H}_{2}(q,Q) &=\frac{1}{2}(Q_{1}^{2}+Q_{2}^{2}+Q_{3}^{2}+Q_{4}^{2})+\frac
{1}{2}(q_{1}^{2}+q_{2}^{2}+q_{3}^{2}+q_{4}^{2}) \; , \\
\Xi(q,Q)&= q_{1}Q_{2}-Q_{1}q_{2} + q_{3}Q_{4}-Q_{3}q_{4}\; , \\
L_{1}(q,Q)&=q_{3}Q_{4}-Q_{3}q_{4} - q_{1}Q_{2}+Q_{1}q_{2}\; ,\\
K(q,Q) &= \frac{1}{2}(-(q_1^2 +Q_1^2)-(q_2^2 +Q_2^2)+(q_3^2 +Q_3^2)+(q_4^2 +Q_4^2)) \; .
\end{align*}
Furthermore
\begin{align*}
N (q,Q)&= \mbox{$\frac{\scriptstyle 1}{\scriptstyle 2}\,$}(K_2^2 +K_3^2) -\mbox{$\frac{\scriptstyle 1}{\scriptstyle 2}\,$}(L_2^2 +L_3^2) \; ,\\
S(q,Q)&=K_{2}L_{3}-K_{3}L_{2} \; ,
\end{align*}
with
\begin{align*}
K_2(q,Q) &=(Q_2Q_3+q_2q_3) -(Q_1Q_4+q_1q_4) \; ,\\
K_3(q,Q) &=-(Q_1Q_3+q_1q_3)-(Q_2Q_4+q_1q_4) \; ,\\
L_2(q,Q) &=(q_1Q_3-Q_1q_3)+(q_2Q_4-Q_2q_4) \; ,\\
L_3(q,Q) &= (q_2Q_3-Q_2q_3)-(q_1Q_4-Q_1q_4) \; .
\end{align*}
This Hamiltonian system has commuting integrals $H_2(q,Q)$, $\Xi(q,Q)$, and $L_1(q,Q)$. That is, $G=\mathbb{T}^3$ generated by these three integrals.
For $G$ we have the orbit map
\[
\rho: (q,Q)\rightarrow (H_2(q,Q),\Xi(q,Q),L_1(q,Q),N(q,Q),K(q,Q),S(q,Q)) \; .
\]

Setting $H_2(q,Q)=n$, $\Xi(q,Q)=\xi$, $L_1(q,Q)=\ell$ we obtain, after reduction with respect to the $\mathbb{T}^3$-action generated by $H_2(q,Q)$, $\Xi(q,Q)$, and $L_1(q,Q)$, the reduced phase space
\[
(n^{2}+\xi^{2}-\ell^{2}-K^{2})^{2}-4(n\xi-\ell K)^{2}=4N^{2}+4S^{2}\; .
\]
in $(N,K,S)$-space. The $\mathbb{T}^3$ momentum map is
\[
{\cal J}: \R^8 \rightarrow (\Xi,L_1,H_2) \subset \R^3 \; ,
\]
which is dual to the orbit map $\rho$. Hence we may classify the symplectic leaves in the orbit space by the values of the momentum map ${\cal J}$, see fig. \ref{RedPhaseSp}.

\begin{figure}
\begin{center}
\includegraphics[width=0.7\columnwidth]{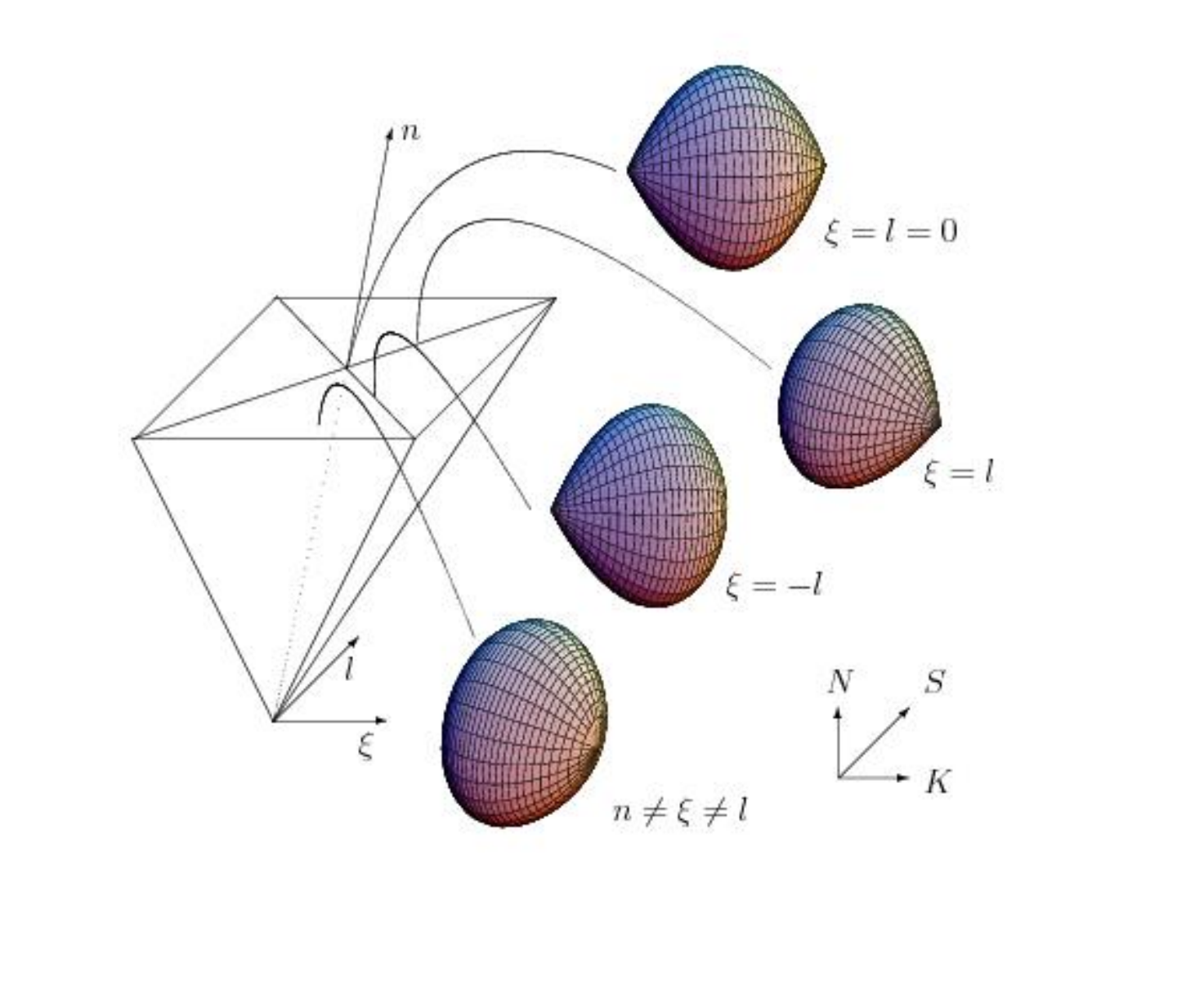}
\caption{Different reduced phase space for the values of ${\cal J}$, \cite{egea2011}}.
\label{RedPhaseSp}
\end{center}
\end{figure}

Note that the image of this momentum map is not a polytope. It is related to a momentum map of deficiency 1 (see \cite{karshon1996}). The image is an upside down pyramid with its diagonal planes

Let $G_{\langle F_1, \cdots ,F_k\rangle}$ denote the group generated by $F_1, \cdots ,F_k$ Consider the action of $F_1=\frac{1}{2} (L_1+\Xi)$. $G_{\langle F_1\rangle}$ is a subgroup of $G_{ \langle H_2, \Xi,L_1\rangle}$, and on $\R^8$ $Fix(G_{\langle F_1\rangle})=\{(q,Q)\in \R^8 |q_1=Q_1=q_2=Q_2=0\}$ is an invariant space. Similarly for the action of $F_2=\frac{1}{2}(\Xi -L_1)$, $G_{\langle F_2\rangle}$ is a subgroup of $G_{\langle H_2, \Xi,L_1\rangle}$, and $Fix(G_{\langle F_2\rangle})=\{(q,Q)\in \R^8 |q_3=Q_3=q_4=Q_4=0\}$ is an invariant space.

${\cal J} (Fix(G_{\langle F_1\rangle}))$ is the restriction of the image of ${\cal J}$ to the plane $\Xi=L_1$.
${\cal J} (Fix(G_{\langle F_2\rangle}))$ is the restriction of the image of ${\cal J}$ to the plane $\Xi=-L_1$.
In the image of the momentum map ${\cal J}$ the fibration in each diagonal plane is equivalent to the fibration of the energy-moment map for the harmonic oscillator. Points in the interior correspond to a fibre topologically equivalent to $T^2$, points on the edges correspond to a fibre topologically equivalent ${S}^1$. A line with $H_2=n$ corresponds to an invariant surface topologically equivalent to ${S}^3$.

The points in the interior of the diagonal planes correspond to the singular reduced phase spaces. The fixed point spaces correspond to the isotropy type and orbit type strata on the reduced phase space which are the zero dimensional symplectic leaves of the final orbit space, that is, they are the cone-like singularities in the singular reduced phase spaces.

The image ${\cal J} (Fix(G_{\langle \Xi\rangle}))$ is given by the planes $\Xi=\pm H_2$.
And ${\cal J}(Fix(G_{\langle L_1\rangle}))$ is given by the planes $L_1=\pm H_2$.
Points in the interior correspond to a fibre topologically equivalent to $T^2$, points on the edges correspond to a fibre topologically equivalent ${S}^1$.

These points correspond to the isotropy strata given by the zero dimensional symplectic leaves on the final orbit space that correspond to the cases where the reduced phase space reduces to an isolated point.

On these singular fibres one finds relative equilibria for all $G$-invariant systems.

Stationary points for the reduced system $X_H$ other then the ones found sofar correspond to the points where the Hamiltonian is tangent to the reduced phase space. These points have as pre-image a $T^3$ on which one finds relative equilibria.

\subsection{$SO(3)$ action on $\R^3$}
This is a classical example. Identify the Lie algebra $so(3)$ of $3\times 3$ matrices with $\R^3$ through
\[
\R^3 \rightarrow so(3);x\rightarrow A_x\; , \;\text{with}\; A_x=\begin{pmatrix}
0 & -x_3 & x_2 \\
x_3 & 0 & -x_1\\
-x_2 & x_1 & 0 \end{pmatrix} \; .
\]
Now $A_x\, y=x\times y$. The Lie bracket on $\R^3$ becomes the cross product, i.e. if $x,y\in \R^3$ then $y\times x \rightarrow [A_x,A_y]=A_xA_y-A_y-A_x$.
Identifying $so(3)$ with its dual $so(3)^*$ we again obtain $\R^3$. The Poisson structure on $so(3)^*$ is $\{ f(x),g(x)\}=\langle x, [Df(x),Dg(x)]\rangle =\langle x,Dg(x)\times Df(x)\rangle $, where $\langle  \; , \; \rangle $ is the standard Euclidian inner product, which is the pairing used to identify $so(3)$ and $so(3)^*$, and $[\; , \;]$ is the Lie bracket on $so(3)$ (cf. the rigid body bracket \cite{marsden1999}). As $\{ f(x),g(x)\}=\langle x,Dg(x)\times Df(x)\rangle =\langle Df(x),x \times Dg(x)\rangle $ the structure matrix for this Poisson structure is $A_x$, which gives also the Poisson structure on $\R^3$ when $so(3)^*$ and $\R^3$ are identified.
The common $SO(3)$ action by $3\times 3$ matrices on $\R^3$ then is a representation of the co-adjoint action of $SO(3)$ on the dual of its Lie algebra. The co-adjoint orbits in $\R^3$ are the spheres $|x|^2=c$, $c\geqslant 0$. Here $|x|^2$ is the Casimir for the Poisson structure on $\R^3$. The only $SO(3)$ invariant is $|x|^2$, which means that the orbit map is $\rho :\R^3 \rightarrow \R ;x\rightarrow |x|^2$. Thus the orbit space is a half-line, which can be considered as a momentum polytope \cite{guillemin1984}.

As $\R^3$ is identified with $so(3)^*$ the momentum map for the $SO(3)$ action on $\R^3$ is the identity $\R^3 \rightarrow \R^3$. And we have a Lie-Weinstein dual pair \cite{weinstein1983,ortega2004} $x \leftarrow x \rightarrow |x|^2$. Moreover, any smooth $f(x)$, $x\in \R^3$ commutes with any smooth $g(|x|^2)$.

Every point $x$ in $\R^3\backslash {0}$ has as isotropy subgroup an $S^1$ given by the rotation around the axis through the origin and $x$. The origin has isotropy subgroup $SO(3)$. The orbit type stratification of the orbit space is given by the origin and the half-line minus the origin.
Every reduced phase space is a point. Every $SO(3)$ invariant Hamiltonian is $f(|x|^2)$ for $f\in C^\infty (\R ,\R )$. Thus every $SO(3)$ invariant Hamiltonian system is integrable and each orbit of an $SO(3)$ invariant Hamiltonian system is a relative equilibrium. A co-adjoint orbit on $\R^3$ reduces to a single point.

\subsection{$SO(3)$ action on $\R^6$}
Consider the action of $SO(3)\times SO(3)$ on $\R^6$ given by $(A,B)\cdot(x,y)=(Ax,By)$. The Lie bracket on $so(3)\times so(3)$ is given by $\{ (x,y),(u,v)\} =(x\times u, y\times v)$. The Poisson bracket on $so(3)^*\times so(3)^*$ after identification with $\R^6$ has structure matrix $\begin{pmatrix} A_x & 0 \\ 0 & A_y \end{pmatrix}$. This Poisson structure is degenerate, the structure matrix has rank 4.
Now restrict the action to its diagonal, that is, consider the $SO(3)$ action on $\R^6$ given by $A\cdot (x,y)=(Ax,Ay)$, with $A\in SO(3)$, $x,y \in \R^3$. Note that this representation is different from the symplectic $SO(3)$ action defined in example \ref{symplso3}.  To put this in a symplectic context consider the co-adjoint orbits in $\R^6$ of the $SO(3)\times SO(3)$ action , which are the symplectic leaves for this Poisson structure. The co-adjoint orbits in $\R^6$ are $S^2\times S^2$ given by $|x|^2=c_x$ and $|y|^2=c_y$.

The invariants for this diagonal $SO(3)$ action are given by $|x|$, $|y|$, the inner product $\langle x,y\rangle $, and the length of the cross product $|x\times y|$. Thus we have the orbit map
\[
\rho :(x,y) \rightarrow (|x|^2, |y|^2, \langle x,y\rangle , |x\times y|)\; .
\]
We have as a relation the Lagrange identity
\begin{equation}\label{lagrangeid}
|x|^2 |y|^2 -\langle x,y\rangle ^2=|x\times y|^2 \; .
\end{equation}

The momentum map for this $SO(3)$-action is
\[
\mu :\R^6 \rightarrow \R^3; (x,y)\rightarrow (\mu_1(x,y),\mu_2(x,y),\mu_3(x,y))=x+ y\; ,
\]
and the components of $x + y$ generate the action of $SO(3)$ on $\R^6$. We have a Lie-Weinstein dual pair $(\mu ,\rho )$
\[
x+ y \overset{\mu}{\longleftarrow} (x,y) \overset{\rho}{\longrightarrow} (|x|^2, |y|^2, \langle x,y\rangle , |x\times y|)\; .
\]
We may consider $x$ and $y$ as two vectors in $\R^3$, which we may interpret as angular momentum vectors. Then $x+y$ is the total angular momentum vector. Only if these two vectors point in the same or in opposite direction we obtain a point $(x,y)$ with isotropy subgroup $S^1$. In this case we have $x=a y$ and $x\times y=0$. The orbit space is the 3-dimensional space in $\R^4$ given by the relation (\ref{lagrangeid}). We may represent this  as solid cone in $\R^3$ with coordinates given by $|x|^2\geqslant 0$, $|y|^2\geqslant 0$, $\langle x,y\rangle $. For each $|x\times y|=c$, we obtain one blade of a two bladed hyperboloid for $c\neq 0$, and a cone if $c=0$.

The Poisson bracket on the image of the orbit map is the induced bracket. A simple computation shows that the brackets between all the invariants vanish, that is, the induced bracket is totally degenerate.

The orbit type stratification has as strata the origin with isotropy subgroup $SO(3)$, the surface of the solid  cone with isotropy subgroup $S^1$ and the interior of the cone with the identity as isotropy subgroup. In this solid cone representation the co-adjoint orbits in $\R^6$ reduce to straight line segments inside the cone, parallel to the $\langle x,y\rangle $-axis, where the end points belong to the cone. When $x=0$ or $y=0$ the co-adjoint orbit is $S^2\times 0$ which reduces to a single point as in the previous section. If we assume conservation of total angular momentum then we obtain the reduced phase spaces. If
\begin{equation}\label{totangmomeq}
|x+y|^2=|x|^2+2\langle x,y\rangle +|y|^2=c_t \; ,
\end{equation}
then we obtain one more equation besides (\ref{lagrangeid}). These two equations define the two dimensional reduced phase space in $R^4$. In the solid cone model this is the intersection of the solid cone with the plane given by (\ref{totangmomeq}). When $c_t >  0$ then this intersection is a parabolic disk. When $c_t=0$ it is a v-shaped plane with the origin in the vertex. When we restrict to the symplectic leaves, that is, the co-adjoint orbits, in $\R^6$, then the values of $|x|^2$ and $|y|^2$ are fixed and, as a consequence of (\ref{totangmomeq}), also $\langle x,y\rangle$ is fixed and the reduced phase space becomes a point.

When we consider the induced Poisson bracket on the image of the moment map, then the structure matrix is $A_\mu=A_x+A_y$, and the co-adjoint orbits on the image of the moment map are $\mu_1^2 +\mu_2^2 +\mu_3^2=c_\mu=c_t$. The pre-image of a co-adjoint orbit on the image of the momentum map in $\R^6$ is then also the pre-image of a reduce phase space according to the one-to-one correspondence between symplectic leaves on the two components of the dual pair \cite{ortega2004}. For each $S^2\times S^2$ the reduced phase space is a line segment and the pre-image of a point on this segment is an $S^2\times S^1$ except for the endpoints which have as pre-image an $S^1\times S^1$.

As an example we may consider the regularized Kepler system \cite{vandermeer2015}(see also subsection \ref{harmosc}). After Kustaanheimo Stiefel regularization the bounded orbits for the Kepler system are in correspondence with the orbits of the harmonic oscillator on $\R^8$ with Hamiltonian
\[
{H}_{2}(q,Q) =\frac{1}{2}(Q_{1}^{2}+Q_{2}^{2}+Q_{3}^{2}+Q_{4}^{2})+\frac
{1}{2}(q_{1}^{2}+q_{2}^{2}+q_{3}^{2}+q_{4}^{2}) \; ,
\]
and integral
\[
\Xi(q,Q)= q_{1}Q_{2}-Q_{1}q_{2} + q_{3}Q_{4}-Q_{3}q_{4}\; .
\]
The Kustaanheimo-Stiefel transformation on the set given by $H_2(q,Q)=1$ and $\Xi (q,Q)=0$ transforms the  momentum integral into $L$ and the eccentricity integral into $K$ given by
\begin{align*}
L_{1}(q,Q)&=q_{3}Q_{4}-Q_{3}q_{4} - q_{1}Q_{2}+Q_{1}q_{2}\; ,\\
L_2(q,Q) &=(q_1Q_3-Q_1q_3)+(q_2Q_4-Q_2q_4) \; ,\\
L_3(q,Q) &= (q_2Q_3-Q_2q_3)-(q_1Q_4-Q_1q_4) \; ,\\
K_1(q,Q) &= \frac{1}{2}(-(q_1^2 +Q_1^2)-(q_2^2 +Q_2^2)+(q_3^2 +Q_3^2)+(q_4^2 +Q_4^2)) \; , \\
K_2(q,Q) &=(Q_2Q_3+q_2q_3) -(Q_1Q_4+q_1q_4) \; ,\\
K_3(q,Q) &=-(Q_1Q_3+q_1q_3)-(Q_2Q_4+q_1q_4) \; .
\end{align*}
These functions generate the $SO(4)$ symmetry of the Kepler system. Define $\sigma _i=L_i+K_1$ and $\rho_i=L_i-K_i$. If we reduce with respect to the $S^1$-actions generated by $H_2$, i.e. the $S^1$ action of the Kepler flow on the bounded orbits, and the $\Xi$-action, the reduced phase space is $S^2 \times S^2$ given by $|\sigma|^2=1$ and $|\rho|^2=1$ (see for instance \cite{moser1970}). We may consider the orbit map $\R^8 \rightarrow \R^6;(q,Q)\rightarrow (\sigma , \rho)$. The induced Poisson structure on $\R^6$ has structure matrix $\begin{pmatrix} A_{\sigma} & 0 \\ 0 & A_{\rho} \end{pmatrix}$ and $\sigma$ and $\rho$ generate an $SO(3)\times SO(3) \cong SO(4)$-action on $\R^6$ with co-adjoint orbit the reduced phase space  $S^2\times S^2$. We may now consider the $SO(3)$  action on $\R^6$ given by the diagonal of $SO(3)\times SO(3)$. The momentum map for this action is precisely $\mu :(q,Q)\rightarrow L(q,Q)$.

\subsection{$SO(3)$ action on $\R^9$}
Next we consider the $SO(3)$ action on $\R^9$ given by $A\cdot (x,y,z)=(Ax,Ay,Az)$, with $A\in SO(3)$, $x,y,z \in \R^3$, being the diagonal of the $SO(3)\times SO(3)\times SO(3)$ action. The Lie bracket on $so(3)\times so(3)\times so(3)$ is given by $\{ (x,y,z),(u,v,w)\} =(x\times u, y\times v,z\times w)$. The Poisson bracket on $so(3)^*\times so(3)^*\times so(3)^*$ after identification with $\R^9$ has structure matrix $\begin{pmatrix} A_x & 0 & 0\\ 0 & A_y & 0\\0&0&A_z \end{pmatrix}$. The co-adjoint orbits in $\R^9$ for the $SO(3)\times SO(3)$  action are $S^2\times S^2\times S^2$ given by $|x|^2=c_x$, $|y|^2=c_y$, and $|z|^2=c_z$. They are symplectic manifolds and invariant under the diagonal of this action. The invariants for the diagonal $SO(3)$ action are given by $|x|$, $|y|$, $|z|$ the inner products $\langle x,y\rangle $, $\langle x,z\rangle $, and $\langle y,z\rangle $, the length of the cross products $|x\times y|$, $|x\times z|$, and $|y \times z|$, and the triple product $\langle (y\times x),z\rangle $. Thus we may consider the orbit map
\begin{multline}
\rho :\R^9 \rightarrow \R^{10};\\(x,y,z) \rightarrow (|x|^2, |y|^2, |z|^2 \langle x,y\rangle , \langle x,z\rangle , \langle y,z\rangle , |x\times y|, |x\times z|, |y\times z|, \langle (y\times x),z\rangle   )\; .
\end{multline}
We have relations
\begin{align}\label{lagrange3}
|x|^2 |y|^2 -\langle x,y\rangle ^2 &=|x\times y|^2 \; , \notag \\
|x|^2 |z|^2 -\langle x,z\rangle ^2 &=|x\times z|^2 \; , \notag \\
|y|^2 |z|^2 -\langle y,z\rangle ^2 &=|y\times z|^2 \; ,
\end{align}
and, using the Gram determinant for the square of the triple product,
\begin{multline}\label{tripleprod}
\langle (y\times x),z\rangle ^2=|x|^2 |y|^2 |z|^2 +2 \langle x,y\rangle \langle y,z\rangle \langle x,z\rangle \\
 -\langle x,z\rangle ^2|y|^2-\langle y,z\rangle ^2|x|^2 -\langle x,y\rangle ^2 |z|^2 \; .
\end{multline}
Thus the image of the orbit map is a six dimensional algebraic set.

We have the moment map
\begin{equation}\mu :\R^9\rightarrow \R^3; (x,y,z)\rightarrow x+y+z \; .
\end{equation}
Because of the relations (\ref{lagrange3}) we may write the orbit map as
\begin{equation}
\rho :\R^9\rightarrow \R^7; (x,y,z) \rightarrow (|x|^2, |y|^2, |z|^2 \langle x,y\rangle , \langle x,z\rangle , \langle y,z\rangle , \langle (y\times x),z\rangle   )\; ,
\end{equation}
where the image is determined as a six dimensional set by (\ref{tripleprod}).
To obtain the reduced phase space we also have to consider total angular momentum constant, that is
\begin{equation}\label{totangmomeq2}
|x|^2+|y|^2 +|z|^2 +2\langle x,y\rangle  +2\langle x,z\rangle +2\langle y,z\rangle =c_t \; ,
\end{equation}
defining a hyperplane intersecting the image of the orbit map giving, if regularity is assumed, a six dimensional reduced phase space.

A symplectic reduction is obtained by restricting to the co-adjoint orbit for the $SO(3)$-action on $\R^9$ which is the manifold  $S^2\times S^2 \times S^2$ given by $|x|^2=c_x$, $|y|^2=c_y$, $|z|^2=c_z$. We may simplify the situation bij assuming the momenta to be non-zero and scaling the vectors $x,y,z$ to length one.

Write $v=(v_1,v_2,v_3,v_4)=(\langle x,y\rangle , \langle x,z\rangle , \langle y,z\rangle , \langle (y\times x),z\rangle   )$. Then the restriction of the orbit map on is $ \rho_s :\R^9\rightarrow \R^4; (x,y,z) \rightarrow (v_1,v_2,v_3,v_4)$. The orbit space is a three dimensional set in $\R^4$ given by
\begin{equation}
v_4^2=1 + 2v_1v_2v_3 -v_1^2-v_2^2 -v_3^2 \; .
 \end{equation}
The reduced phase spaces are obtained by intersecting the orbit space with the total momentum plane given by
\begin{equation}\label{totangmomeq3}
\frac{1}{\sqrt{c_xc_y}}v_1 +\frac{1}{\sqrt{c_xc_z}}v_2+\frac{1}{\sqrt{c_yc_z}}v_3=C \;,
\end{equation}
and have maximally dimension two.
In this case the structure matrix for the induced Poisson bracket on the orbit space is given by

\begin{multline}\label{StructMatrR4}
W(v)=\\
\left(
\begin{array}{cccc}
0 & \frac{v_4}{c_x} & -\frac{v_4}{c_y} & \frac{v_1v_3-v_2}{c_x}-\frac{v_1v_2-v_3}{c_y}  \\
-\frac{v_4}{c_x} & 0 &  \frac{v_4}{c_z}&\frac{v_1-v_2v_3}{c_x}+\frac{v_1v_2-v_3}{c_z}\\
\frac{v_4}{c_y} & -\frac{v_4}{c_z} & 0 &\frac{v_2v_3+v_1}{c_y} -\frac{v_1v_3-v_2}{c_z}\\
-\frac{v_1v_3-v_2}{c_x}+\frac{v_1v_2-v_3}{c_y}  &-\frac{v_1-v_2v_3}{c_x}-\frac{v_1v_2-v_3}{c_z} & -\frac{v_2v_3+v_1}{c_y} +\frac{v_1v_3-v_2}{c_z} & 0 \\
\end{array}
\right) \; .
\end{multline}

A system with this symmetry is plays a role in the averaged full two body problem \cite{boue2009} and in compact planetary systems perturbed by an inclined companion \cite{boue2014}. By averaging a Hamiltonian system is obtained in the invariants $\langle x,y\rangle , \langle x,z\rangle , \langle y,z\rangle$. When these are taken as the new coordinates $(v_1,v_2,v_3)$ then, as one can see from (\ref{StructMatrR4}), a Hamiltonian system is obtained  that has a Poisson structure that is totally degenerate along $|\langle (y\times x),z\rangle|^2=0$, which is also the condition for the vectors $x,y,z$ to be coplanar. If we restrict to a co-adjoint orbit in $\R^9$ the surface given by $v_4=0$ is in $(v_1,v_2,v_3)$ space a sort of inflated tetrahedron, part of a Cayley’s nodal cubic surface, also called an elliptope or a Cassini berlingot shaped volume. It is smooth except at the four vertices of the tetrahedron which correspond to the four collinear configurations of the vectors $x,y,z$. Taking into account the conservation of total angular momentum given by (\ref{totangmomeq2}) a two dimensional planar phase space is obtained bounded by a curve, which is the intersection of the total angular momentum plane and the elliptope, and along which the symplectic form vanishes. We have a folded symplectic or b-symplectic form \cite{guillemin2014,dasilva2000}. However, in approach of \cite{boue2009}, \cite{boue2014} the invariant $v_4$ is neglected. If the symmetry is recognized and the full set of invariants is taken into account the dynamics  becomes much easier and a full reduction can be performed, in which case the planar phase space is unfolded to a two dimensional sphere-like, possibly non-smooth, reduced phase space on which the dynamics easily can be constructed.

If we restrict to a co-adjoint orbit in $\R^9$ we may use equation (\ref{totangmomeq3}) to eliminate $v_3$, to obtain a visualisation of the reduced phase space in $(v_1,v_2,v_4)$-space. An example is given in figure \ref{fig:RedSyst3}.

\begin{figure}[htb]
\begin{centering}
\includegraphics[scale=0.52]{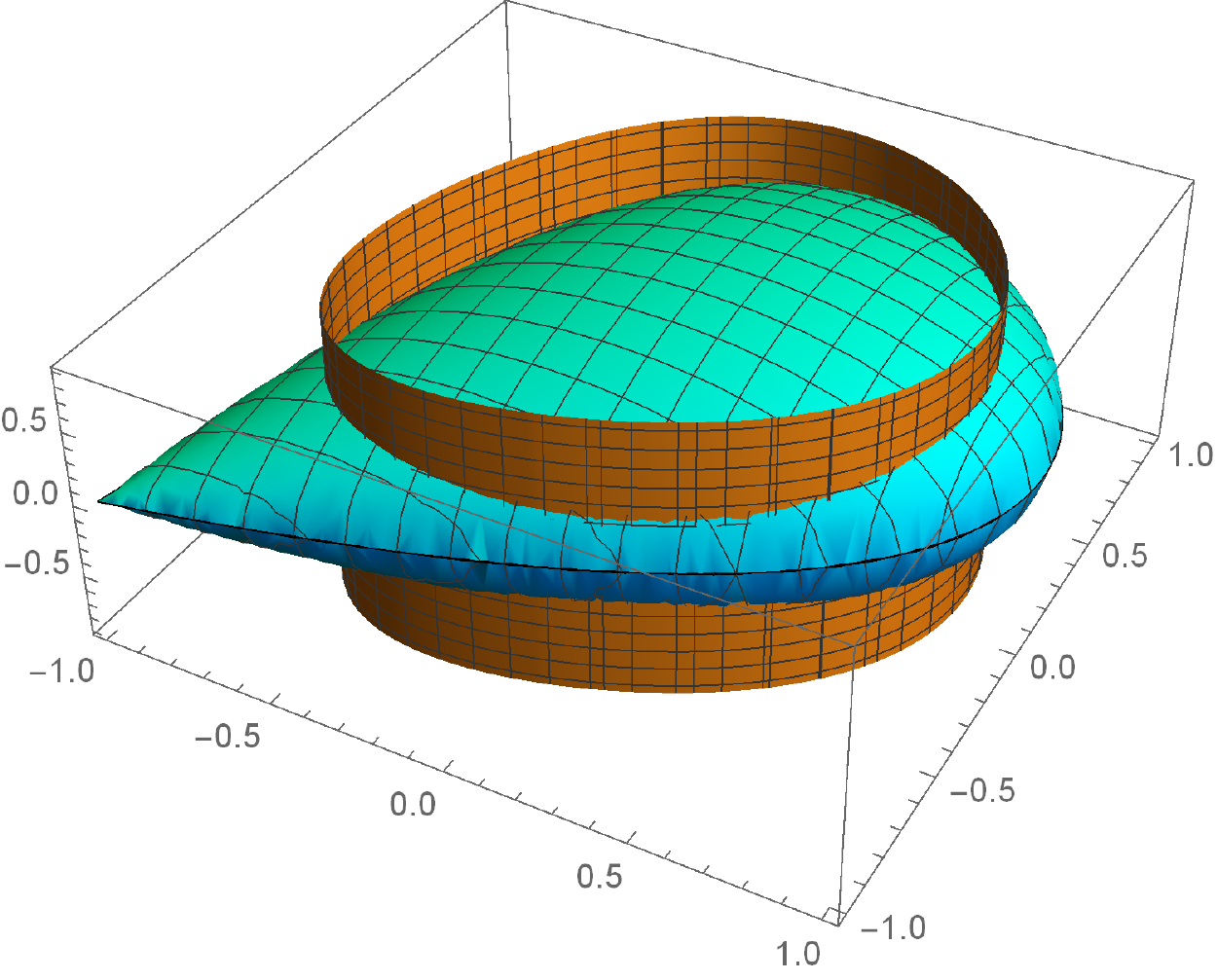}
\caption{Example of non-smooth reduced phase space.}
\label{fig:RedSyst3}
\end{centering}
\end{figure}

The reduced phase space can be considered as the set of possible configurations of the three vectors $x,y,z$. The reduced phase space has a singular point if the curve $v_4=0$ contains a collinear situation. In figure \ref{fig:RedSyst3}. the cylindrical surface parallel to the vertical $v_4$-axis represents a level surface of a reduced Hamiltonian $H(v_1,v_2)$ which does not depend on $v_4$. The intersections of the level surface and the reduced phase space are the trajectories of the reduced system. The reduced Hamiltonian system $X_H$ has stationary points when these level surfaces are tangent to the reduced phase space or when the level surface becomes a vertical line, that is when $D_vH(v)=0$. Stationary points of the reduced system have as pre-image a group orbit which is filled with $X_H$ relative equilibria. The pre-image of a singular point on the reduced phase space has, because of the collinear situation an $S^1$ isotropy subgroup. The singular point belongs to a different isotropy type stratum. The problem of the full averaged model of two interacting rigid bodies will be studied in detail in the context of reduction by invariants in \cite{crespo2023} including a further analysis of the above reduction.

\end{document}